\documentclass[conference]{IEEEtran}
\IEEEoverridecommandlockouts
\usepackage{cite}
\usepackage{amsmath,amssymb,amsfonts,amsthm}
\ifCLASSOPTIONcompsoc
    \usepackage[caption=false, font=normalsize, labelfont=sf, textfont=sf]{subfig}
\else
\usepackage[caption=false, font=footnotesize]{subfig}
\fi
\usepackage[dvipsnames]{xcolor}
\definecolor{mycolor1}{rgb}{1, 0, 0}
\definecolor{mycolor2}{rgb}{0, 0, 1}

\usepackage{tikz}
\usepackage{pgfplots}
\usepackage{algorithmic}
\usepackage[linesnumbered,ruled,vlined]{algorithm2e}
\SetInd{0.25em}{0.25em}
\usepackage{graphicx}
\usepackage{textcomp}
\usepackage{hyperref}
\def\BibTeX{{\rm B\kern-.05em{\sc i\kern-.025em b}\kern-.08em
    T\kern-.1667em\lower.7ex\hbox{E}\kern-.125emX}}

\newtheorem{proposition}{Proposition}
\newtheorem{assumption}{Assumption}

\newcommand\mydots{\makebox[1em][c]{.\hfil.\hfil.}}

\newcommand{\funcblock}[1]{\ensuremath{f_{#1}}}
\newcommand{\funcaltblock}[1]{\ensuremath{g_{#1}}}
\newcommand{\minimizer}{\textcolor{mycolor2}{\ensuremath{a}}}
\newcommand{\maximizer}{\textcolor{mycolor1}{\ensuremath{b}}}
\newcommand{\dims}{\ensuremath{n}}
\newcommand{\dimminimizer}{\textcolor{mycolor2}{\ensuremath{\dims_{\minimizer}}}}
\newcommand{\dimmaximizer}{\textcolor{mycolor1}{\ensuremath{\dims_{\maximizer}}}}
\newcommand{\dimminimizerblock}[1]{\textcolor{mycolor2}{\ensuremath{\dims_{\minimizer,#1}}}}
\newcommand{\dimmaximizerblock}[1]{\textcolor{mycolor1}{\ensuremath{\dims_{\maximizer,#1}}}}
\newcommand{\dualvar}{\ensuremath{\lambda}}
\newcommand{\varminimizer}{\textcolor{mycolor2}{\ensuremath{x_{\minimizer}}}}
\newcommand{\varmaximizer}{\textcolor{mycolor1}{\ensuremath{x_{\maximizer}}}}
\newcommand{\auxvarminimizer}{\textcolor{mycolor2}{\ensuremath{z_{\minimizer}}}}
\newcommand{\auxvarmaximizer}{\textcolor{mycolor1}{\ensuremath{z_{\maximizer}}}}
\newcommand{\dualvarminimizer}{\textcolor{mycolor2}{\ensuremath{\dualvar_{\minimizer}}}}
\newcommand{\dualvarmaximizer}{\textcolor{mycolor1}{\ensuremath{\dualvar_{\maximizer}}}}
\newcommand{\spaceminimizer}{\textcolor{mycolor2}{\ensuremath{\mathcal{X}_{\minimizer}}}}
\newcommand{\spacemaximizer}{\textcolor{mycolor1}{\ensuremath{\mathcal{X}_{\maximizer}}}}

\newcommand{\varminimizerblock}[1]{\textcolor{mycolor2}{\ensuremath{x_{\minimizer, #1}}}}
\newcommand{\varmaximizerblock}[1]{\textcolor{mycolor1}{\ensuremath{x_{\maximizer, #1}}}}
\newcommand{\auxvarminimizerblock}[1]{\textcolor{mycolor2}{\ensuremath{z_{\minimizer, #1}}}}
\newcommand{\auxvarmaximizerblock}[1]{\textcolor{mycolor1}{\ensuremath{z_{\maximizer, #1}}}}

\newcommand{\spaceminimizerblock}[1]{\textcolor{mycolor2}{\ensuremath{\mathcal{X}_{\minimizer, #1}}}}
\newcommand{\spacemaximizerblock}[1]{\textcolor{mycolor1}{\ensuremath{\mathcal{X}_{\maximizer, #1}}}}

\newcommand{\dualstepsizeminimizer}{\textcolor{mycolor2}{\rho_{\minimizer}}}
\newcommand{\dualstepsizemaximizer}{\textcolor{mycolor1}{\rho_{\maximizer}}}
\newcommand{\ltwonorm}[1]{\left\lVert #1 \right\rVert^{2}_{2}}
\newcommand{\ltwonormnormal}[1]{\left\lVert #1 \right\rVert_{2}}
\newcommand{\indicatorfunc}[1]{I_{#1}}
\newcommand{\optvalue}{p^{*}}

\newcommand{\varminimizeritr}[1]{\textcolor{mycolor2}{\ensuremath{x_{\minimizer}^{#1}}}}
\newcommand{\varmaximizeritr}[1]{\textcolor{mycolor1}{\ensuremath{x_{\maximizer}^{#1}}}}
\newcommand{\auxvarminimizeritr}[1]{\textcolor{mycolor2}{\ensuremath{z_{\minimizer}^{#1}}}}
\newcommand{\auxvarmaximizeritr}[1]{\textcolor{mycolor1}{\ensuremath{z_{\maximizer}^{#1}}}}
\newcommand{\dualvarminimizeritr}[1]{\textcolor{mycolor2}{\ensuremath{\dualvar_{\minimizer}^{#1}}}}
\newcommand{\dualvarmaximizeritr}[1]{\textcolor{mycolor1}{\ensuremath{\dualvar_{\maximizer}^{#1}}}}
\newcommand{\primalresidualminimizeritr}[1]{\textcolor{mycolor2}{\ensuremath{r_{\minimizer}^{#1}}}}
\newcommand{\primalresidualmaximizeritr}[1]{\textcolor{mycolor1}{\ensuremath{r_{\maximizer}^{#1}}}}

\newcommand{\varminimizerblockitr}[2]{\textcolor{mycolor2}{\ensuremath{x_{\minimizer, #1}^{#2}}}}
\newcommand{\varmaximizerblockitr}[2]{\textcolor{mycolor1}{\ensuremath{x_{\maximizer, #1}^{#2}}}}
\newcommand{\auxvarminimizerblockitr}[2]{\textcolor{mycolor2}{\ensuremath{z_{\minimizer, #1}^{#2}}}}
\newcommand{\auxvarmaximizerblockitr}[2]{\textcolor{mycolor1}{\ensuremath{z_{\maximizer, #1}^{#2}}}}
\newcommand{\dualvarminimizerblockitr}[2]{\textcolor{mycolor2}{\ensuremath{\dualvar_{\minimizer, #1}^{#2}}}}
\newcommand{\dualvarmaximizerblockitr}[2]{\textcolor{mycolor1}{\ensuremath{\dualvar_{\maximizer, #1}^{#2}}}}

\newcommand{\valueitr}[1]{p^{#1}}
\newcommand{\optvalueminimizeritr}[1]{\textcolor{mycolor1}{(p^{*}_{\maximizer})^{#1}}}
\newcommand{\optvaluemaximizeritr}[1]{\textcolor{mycolor2}{(p^{*}_{\minimizer})^{#1}}}
\newcommand{\valuefunctionitr}[1]{V^{#1}}
\begin{document}

\title{Alternating Direction Method of Multipliers for Decomposable Saddle-Point Problems\\
\thanks{This work was supported in part by NSF 1652113 and ARO W911NF-20-1-0140.}
}

\author{\IEEEauthorblockN{Mustafa O. Karabag}
\IEEEauthorblockA{\textit{Electrical \& Computer Engineering} \\
\textit{The University of Texas at Austin}\\
Austin, USA \\
karabag@utexas.edu}
\and
\IEEEauthorblockN{David Fridovich-Keil}
\IEEEauthorblockA{\textit{
Aerospace Engineering} \\
\textit{The University of Texas at Austin}\\
Austin, USA \\
dfk@utexas.edu}
\and
\IEEEauthorblockN{Ufuk Topcu}
\IEEEauthorblockA{\textit{Aerospace Engineering} \\
\textit{The University of Texas at Austin}\\
Austin, USA \\
utopcu@utexas.edu}
}

\maketitle

\begin{abstract}
Saddle-point problems appear in various settings including machine learning, zero-sum stochastic games, and regression problems. We consider decomposable saddle-point problems and study an extension of the alternating direction method of multipliers to such saddle-point problems. Instead of solving the original saddle-point problem directly, this algorithm solves smaller saddle-point problems by exploiting the decomposable structure. We show the convergence of this algorithm for convex-concave saddle-point problems under a mild assumption. We also provide a sufficient condition for which the assumption holds. We demonstrate the convergence properties of the saddle-point alternating direction method of multipliers with numerical examples on a power allocation problem in communication channels and a network routing problem with adversarial costs.
\end{abstract}


\section{Introduction}
Saddle-point problems consider optimization of an objective function simultaneously by a minimizer and maximizer. These problems appear, for example, in zero-sum stochastic games~\cite{shapley1953stochastic}, adversarial training of machine learning models~\cite{goodfellow2014generative,sinha2017certifying}, regression problems~\cite{xu2009robustness}, and maximum-margin estimation of structured output models~\cite{taskar2006structured}.

We focus on decomposable saddle-point problems of the following form that have a decomposable objective function with \textit{complicating} global constraints:
\begin{subequations} \label{prob:original}
\begin{align} 
\min_{\varminimizer} \max_{\varmaximizer} \quad & \sum_{i=1}^{N} \funcblock{i}(\varminimizerblock{i}, \varmaximizerblock{i}) 
\\
\text{subject to} \quad & \varminimizer \in \spaceminimizer \label{cons:minglobal}
\\
&\varmaximizer \in \spacemaximizer \label{cons:maxglobal}
\\
&\varminimizerblock{i} \in \spaceminimizerblock{i},  \quad \text{for all } i \in 1, \ldots, N 
\\
&\varmaximizerblock{i} \in \spacemaximizerblock{i}, \quad  \text{for all } i \in 1, \ldots, N 
\end{align}
\end{subequations} where \(\varminimizer = [\varminimizerblock{1}, \ldots, \varminimizerblock{N}]\) and \(\varmaximizer = [\varmaximizerblock{1}, \ldots, \varmaximizerblock{N}]\) are each concatenations of \(N\) vectors, and \(\spaceminimizer \subseteq \mathbb{R}^{\dimminimizer}\), \(\spacemaximizer \subseteq \mathbb{R}^{\dimmaximizer}\), \(\spaceminimizerblock{i} \subseteq \mathbb{R}^{\dimminimizerblock{i}},\) and \(\spacemaximizerblock{i} \subseteq \mathbb{R}^{\dimmaximizerblock{i}}\) are compact, convex sets such that \(\sum_{i = 1}^{N} \dimminimizerblock{i} = \dimminimizer \) and \(\sum_{i = 1}^{N} \dimmaximizerblock{i} = \dimmaximizer.\) In particular, we are interested in the convex-concave case, i.e., \(\funcblock{i}: \spaceminimizerblock{i}\times \spacemaximizerblock{i} \to \mathbb{R}\) is convex, lower semicontinuous in \(\varminimizerblock{i}\) and concave, upper semicontinuous in \(\varmaximizerblock{i}\). This problem structure arises, for example, in power allocation problems for communication channels with adversarial noise~\cite{ghosh2003minimax} and optimal network routing problems with adversarial costs.

 The paper \cite{karabag2022deception} proposed the alternating direction method of multipliers (ADMM) to solve an optimization problem with decomposable nonconvex-concave objective functions. In this paper, we analyze the convergence properties of this method, saddle-point ADMM (SP-ADMM), for the decomposable convex-concave objective functions. The iterative SP-ADMM preserves the separable structure of \eqref{prob:original} and consists of three steps. In the first step, SP-ADMM solves a saddle-point problem separately for every block.  It performs projections onto the global constraints \eqref{cons:minglobal}--\eqref{cons:maxglobal} in the next step, and performs the dual variable updates in the last step. 
 
 SP-ADMM has several advantages. Each individual saddle-point problem has a lower number of dimensions compared to the original problem and hence can be solved more efficiently. For some objective functions, for example bilinear functions of two one-dimensional variables, these individual saddle-point problems can be solved analytically. Since the individual saddle-point problems have no coupling, they can be solved in parallel. SP-ADMM performs the projection onto the global constraints without considering the individual constraints. For some  global constraints such as unit ball or probability simplex, this projection step can be performed more efficiently compared to the case that takes the individual constraints into account.
 
 The contributions of this paper are threefold. The paper \cite{karabag2022deception} demonstrated the performance of SP-ADMM for a specific robust optimization problem without any theoretical guarantees. We analyze the performance of SP-ADMM. We first show that for the convex-concave case SP-ADMM converges to the saddle point of the problem under a mild assumption. Secondly, we provide a sufficient condition for convergence by considering standard conditions of the minimax theorem~\cite{border1985fixed} and Slater's constraint qualification~\cite{boyd2004convex}. Finally, we demonstrate and evaluate the performance of SP-ADMM for a power allocation problem for communication channels with adversarial noise~\cite{ghosh2003minimax} and an optimal network routing problem with adversarial costs.

\section{Related Work}
\textit{Saddle-point problems:} 
Convergent variants of gradient descent-ascent methods such as the extra gradient method~\cite{korpelevich1976extragradient}, optimistic gradient descent-ascent method~\cite{daskalakis2018limit}, and subgradient descent-ascent method~\cite{nedic2009subgradient} have been proposed for convex-concave saddle-point problems. The paper \cite{gidel2017frank} extended the Frank-Wolfe (conditional gradient) method to solve strongly convex-strongly concave saddle-point problems. 

We remark that first-order methods can also exploit the decomposable structure during gradient computation. However, the coexistence of local and global constraints for the projection step may result in harder optimization problems compared to SP-ADMM that decouples the  projection step and local constraints. 

SP-ADMM solves saddle-point problems as a subroutine and one can employ these methods to solve the individual saddle-problems. The quadratic penalties introduced in the SP-ADMM results in strongly convex-strongly concave objective functions that often increase the rate of convergence.

\textit{Decomposable optimization:} 
Decomposable optimization studies optimization problems that can be decomposed into smaller sub-problems once the complicating constraints (or variables) are removed. Seminal Dantzig–Wolfe~\cite{dantzig1960decomposition} and Benders~\cite{bnnobrs1962partitioning} decomposition methods solve block decomposable linear programs. ADMM~\cite{gabay1976dual,glowinski1975approximation} solves general decomposable convex optimization problems. ADMM has convergence guarantees for convex problems~\cite{boyd2011distributed,nishihara2015general} and also for some nonconvex problems~\cite{guo2017convergence,wang2019global}. In practice, ADMM often generates acceptable solutions in a few iterations, however it behaves like a first-order method and suffers from slow convergence in the long run~\cite{boyd2011distributed}.

\textit{Decentralized saddle-point problems:} 
Decentralized saddle-point problems~\cite{liu2019decentralized,hou2021efficient,rogozin2021decentralized,sharma2022federated,mateos2015distributed} consider the optimization of a separable objective function subject to the communication constraints (usually defined with a graph). Unlike the decomposable setting that we consider, these works consider that each component of the objective function is a function of a global variable. The paper \cite{rogozin2021decentralized} also consider local variables as a part of the objective functions, however these local variables do not have complicating constraints that we have in \eqref{prob:original}.
\section{Notation and Preliminaries for Decomposable Optimization}

\subsection{Notation} We use subscripts \(\minimizer\) and \(\maximizer\) with colors blue and red to denote the variables/constants of the minimizer and maximizer, respectively. The subscript \(i\) denotes the \(i^{\text{th}}\) block (element) of the object with the subscript. With an abuse of notation we also use the subscript \(i\), for \(\spaceminimizerblock{i}\) and \(\spacemaximizerblock{i}\) that the sets for \(\varminimizerblock{i}\) and \(\varmaximizerblock{i}\), and are not blocks of \(\spaceminimizer\) and \(\spacemaximizer\), respectively. The superscript \(k\) denotes the value of the variable with the superscript at the \(k^{\text{th}}\) iteration of the algorithms. The superscript \(2\) is used as an exponent. \(\indicatorfunc{\mathcal{X}}(x) \) is the indicator function of set \(\mathcal{X}\) such that \(\indicatorfunc{\mathcal{X}}(x)=0\) if \(x \in \mathcal{X}\) and \(\infty\) otherwise.

\subsection{Preliminaries for Decomposable Optimization} 
While the objective function is block-decomposable for \eqref{prob:original}, the constraints \(\varminimizer \in \spaceminimizer\) and \(\varmaximizer \in \spacemaximizer\) are not separable. The potential existence of these constraints result in different problem structures: 
 \paragraph{Fully separable case} In the absence of both \(\varminimizer \in \spaceminimizer\) and \(\varmaximizer \in \spacemaximizer\), we can solve the saddle-point problem separately for every block \(i\). The saddle points of these individual problems are jointly a saddle point for the global problem.
     \paragraph{Maximizer separable case} In this case, the inner maximization problem is a function of \(\varminimizer\), i.e.,  \(\hat{\funcblock{i}}(\varminimizerblock{i}) = \max_{\varmaximizerblock{i} \in \spacemaximizerblock{i}} \funcblock{i}(\varminimizerblock{i}, \varmaximizerblock{i})\). If \(\hat{\funcblock{i}}\) can be derived, we get a minimization problem with a block-separable objective. However, the global minimization problem still contains constraint \(\varmaximizer \in \spacemaximizer\). This optimization problem can be solved with decomposable optimization methods such as ADMM.
    \paragraph{Inseparable case} If both \(\varminimizer \in \spaceminimizer\) and \(\varmaximizer \in \spacemaximizer\) are present, we cannot use  \(\hat{\funcblock{i}}(\varminimizerblock{i}) =\max_{\varmaximizerblock{i} \in \spacemaximizerblock{i}} \funcblock{i}(\varminimizerblock{i}, \varmaximizerblock{i})\) due to the globally bounding constraint \(\varmaximizer \in \spacemaximizer\). We may attempt to derive \(\hat{f}(\varminimizer) =\max_{\varmaximizer \in \spacemaximizer \cap (\cap_{i=1}^{N} \spacemaximizerblock{i})}  \funcblock{i}(\varminimizerblock{i}, \varmaximizerblock{i})\). However, this process (potentially) removes the separability of the objective function. Hence, decomposable optimization methods are not directly applicable to this case. We are interested in separable solutions for this case by preserving the minimax formulation.

\section{Alternating Direction Method of Multipliers for Decomposable Optimization}

The alternating direction method of multipliers (ADMM)~\cite{gabay1976dual,glowinski1975approximation} is an optimization method to solve optimization problems with separable objectives and complicating constraints. Consider the problem
\begin{subequations} \label{prob:onlyminoriginal}
\begin{align} 
\min_{\varminimizer} \quad & \sum_{i=1}^{N} \funcaltblock{i}(\varminimizerblock{i}) 
\\
\text{subject to} \quad & \varminimizer \in \spaceminimizer
\\
&\varminimizerblock{i} \in \spaceminimizerblock{i},  \quad \text{for all } i \in 1, \ldots, N.
\end{align}
\end{subequations}
To solve this problem using ADMM, we use an auxiliary variable \(\auxvarminimizer\) and rewrite \eqref{prob:onlyminoriginal} as 
\begin{subequations} \label{prob:onlyminwithaux}
\begin{align} 
\min_{\varminimizer} \quad & \sum_{i=1}^{N} \funcaltblock{i}(\varminimizerblock{i}) 
\\
\text{subject to} \quad & \auxvarminimizer \in \spaceminimizer
\\
&\auxvarminimizerblock{i} = \varminimizerblock{i},  \quad \text{for all } i \in 1, \ldots, N
\\
&\varminimizerblock{i} \in \spaceminimizerblock{i},  \quad \text{for all } i \in 1, \ldots, N.
\end{align}
\end{subequations}
For \eqref{prob:onlyminwithaux}, we define the Lagrangian \(\mathcal{K}(\varminimizer, \auxvarminimizer, \dualvarminimizer)\) as
\begin{align*}
      \sum_{i=1}^{N} \left(\funcaltblock{i}(\varminimizerblock{i}) + \indicatorfunc{\spaceminimizerblock{i}}  (\varminimizerblock{i})  \right) + \dualvarminimizer^{\top}(\varminimizer - \auxvarminimizer)  + \indicatorfunc{\spaceminimizer}(\auxvarminimizer) 
\end{align*} and the augmented Lagrangian \(\hat{\mathcal{K}}(\varminimizer, \auxvarminimizer, \dualvarminimizer)\) as 
\begin{align*}
  \mathcal{K}(\varminimizer, \auxvarminimizer, \dualvarminimizer) +  \frac{\dualstepsizeminimizer}{2} \ltwonorm{\varminimizer - \auxvarminimizer}
\end{align*}
where \(\dualstepsizeminimizer > 0 \) is the penalty parameter.

\begin{algorithm}[t]
    \DontPrintSemicolon 
    \SetKwBlock{DoParallel}{For every \(i \in [\numAgents]\) do in parallel}{end}
    
    Initialize \(\varminimizeritr{0}, \auxvarminimizeritr{0}, \dualvarminimizeritr{0}\) such that \(\varminimizeritr{0} \in \spaceminimizerblock{1} \times \ldots \times \spaceminimizerblock{N}\), \(\auxvarminimizeritr{0} \in \spaceminimizer\).

    \For{\( k = 0,1,\ldots \)}{
\(\varminimizeritr{k+1} = \arg \min_{\varminimizer }  \ \hat{\mathcal{K}}(\varminimizer, \auxvarminimizeritr{k},  \dualvarminimizeritr{k}) .\)  \label{step:minADMMvarupdate}
        
    \(\auxvarminimizeritr{k+1} = \arg \min_{\auxvarminimizer} \ \hat{\mathcal{K}}(\varminimizeritr{k+1},  \auxvarminimizer, \dualvarminimizeritr{k}).\) \label{step:minADMMauxvarupdate} 
         
\(\dualvarminimizeritr{k+1} = \dualvarminimizeritr{k} + \dualstepsizeminimizer(\varminimizeritr{k+1} - \auxvarminimizeritr{k+1})\) \label{step:minADMMdualvarupdate}

    }

    \caption{Alternating Direction Method of Multipliers (ADMM) for decomposable optimization}
    \label{alg:admm}
\end{algorithm}

ADMM for decomposable optimization, Algorithm \ref{alg:admm}, consists of three steps: primal variable \(\varminimizer\), auxiliary primal variable \(\auxvarminimizer\), and dual variable \(\dualvarminimizer\) updates. We note that Line \ref{step:minADMMvarupdate} of Algorithm \ref{alg:admm} is separable and is the same with assigning  \[ \arg \min_{\varminimizerblock{i} \in \spaceminimizerblock{i}} \ \funcaltblock{i}(\varminimizerblock{i})  + (\dualvarminimizerblockitr{i}{k})^{\top}(\varminimizerblock{i} - \auxvarminimizerblockitr{i}{k}) + \frac{\dualstepsizeminimizer}{2} \ltwonorm{\varminimizerblock{i} - \auxvarminimizerblockitr{i}{k}}  \] to \(\varminimizerblockitr{i}{k+1}\) for every \(i \in 1, \ldots N\). Line \ref{step:minADMMauxvarupdate} is the convex projection step and is equal to letting \[\auxvarminimizeritr{k+1} = \arg \min_{\auxvarminimizer \in \spaceminimizer}  (\dualvarminimizeritr{k})^{\top}(\varminimizeritr{k+1} - \auxvarminimizer) + \frac{\dualstepsizeminimizer}{2} \ltwonorm{\varminimizeritr{k+1} - \auxvarminimizer} \] which is equal to \[\auxvarminimizeritr{k+1} = \arg \min_{\auxvarminimizer \in \spaceminimizer} \ltwonorm{\varminimizeritr{k+1} + \dualvarminimizeritr{k}/\dualstepsizeminimizer - \auxvarminimizer} .\]

Let \(\varminimizeritr{*}\) be an optimal solution of \eqref{prob:onlyminoriginal}. The iterates of ADMM converge to an optimal solution\cite{boyd2011distributed}, i.e., \(\auxvarminimizeritr{k} \to \varminimizeritr{*}\), if there exists \((\varminimizeritr{*}, \auxvarminimizeritr{*},  \dualvarminimizeritr{*})\) for all \(\varminimizer\), \(\auxvarminimizer\), and \(\dualvarminimizer\) such that 
\begin{equation} \label{ineq:saddlepointmin}
    \mathcal{K} (\varminimizeritr{*}, \auxvarminimizeritr{*}, \dualvarminimizer) \leq \mathcal{K} (\varminimizeritr{*}, \auxvarminimizeritr{*}, \dualvarminimizeritr{*}) \leq     \mathcal{K} (\varminimizer,  \auxvarminimizer,  \dualvarminimizeritr{*}).
\end{equation}

\section{Saddle-Point Alternating Direction Method of Multipliers}
In this section, we describe the alternating direction method of multipliers (ADMM) for saddle-point problems that was first introduced in \cite{karabag2022deception}. The method shares the same steps with standard ADMM and enjoys the same convergence guarantees.

To apply ADMM to saddle-point problem \eqref{prob:original} we first rewrite the problem using the auxiliary variables \(\auxvarminimizerblock{i}\) and \(\auxvarmaximizerblock{i}\):
\begin{subequations} \label{prob:withaux}
\begin{align} 
\min_{\varminimizer, \auxvarminimizer} \max_{\varmaximizer, \auxvarmaximizer}\quad & \sum_{i=1}^{N} \funcblock{i}(\varminimizerblock{i}, \varmaximizerblock{i}) \\ 
\text{subject to} \quad & \auxvarminimizer \in \spaceminimizer
\\
&\auxvarmaximizer \in \spacemaximizer
\\
&\auxvarminimizerblock{i} = \varminimizerblock{i},  \quad \text{for all } i \in 1, \ldots, N,
\\
&\auxvarmaximizerblock{i} = \varmaximizerblock{i}, \quad  \text{for all } i \in 1, \ldots, N,
\\
&\varminimizerblock{i} \in \spaceminimizerblock{i},  \quad \text{for all } i \in 1, \ldots, N,
\\
&\varmaximizerblock{i} \in \spacemaximizerblock{i}, \quad  \text{for all } i \in 1, \ldots, N.
\end{align}
\end{subequations}

For \eqref{prob:withaux}, we define the Lagrangian 
\begin{align*}
    \mathcal{L}&(\varminimizer, \varmaximizer, \auxvarminimizer, \auxvarmaximizer, \dualvarminimizer, \dualvarmaximizer)=
    \\
    &  \sum_{i=1}^{N} \left(\funcblock{i}(\varminimizerblock{i}, \varmaximizerblock{i}) + \indicatorfunc{\spaceminimizerblock{i}}  (\varminimizerblock{i}) - \indicatorfunc{\spacemaximizerblock{i} }(\varmaximizerblock{i}) \right)
    \\
    &+ \dualvarminimizer^{\top}(\varminimizer - \auxvarminimizer)  + \indicatorfunc{\spaceminimizer}(\auxvarminimizer) - \dualvarmaximizer^{\top}(\varmaximizer - \auxvarmaximizer)  - \indicatorfunc{\spacemaximizer}(\auxvarmaximizer) 
\end{align*} and the augmented Lagrangian 
\begin{align*}
    \hat{\mathcal{L}}(\varminimizer, \varmaximizer, \auxvarminimizer, \auxvarmaximizer, \dualvarminimizer, \dualvarmaximizer) & =    \mathcal{L}(\varminimizer, \varmaximizer, \auxvarminimizer, \auxvarmaximizer, \dualvarminimizer, \dualvarmaximizer)
    \\ &+ \frac{\dualstepsizeminimizer}{2} \ltwonorm{\varminimizer - \auxvarminimizer} - \frac{\dualstepsizemaximizer}{2} \ltwonorm{\varmaximizer - \auxvarmaximizer} .
\end{align*}
where \(\dualstepsizeminimizer > 0 \) is the penalty parameter for the minimizer, and \(\dualstepsizemaximizer > 0 \) is the penalty parameter for the maximizer.

\begin{algorithm}[t]
    \DontPrintSemicolon 
    \SetKwBlock{DoParallel}{For every \(i \in [\numAgents]\) do in parallel}{end}
    
    Initialize \(\varminimizeritr{0}, \varmaximizeritr{0}, \auxvarminimizeritr{0}, \auxvarmaximizeritr{0}, \dualvarminimizeritr{0}, \dualvarmaximizeritr{0}\) such that \(\varminimizeritr{0} \in \spaceminimizerblock{1} \times \ldots \times \spaceminimizerblock{N}\), \(\varmaximizeritr{0} \in \spacemaximizerblock{1} \times \ldots \times \spacemaximizerblock{N}\), \(\auxvarminimizeritr{0} \in \spaceminimizer\), and \(\auxvarmaximizeritr{0} \in \spacemaximizer\).

    \For{\( k = 0,1,\ldots \)}{
        
            \( \varminimizeritr{k+1}, \varmaximizeritr{k+1} = \arg \underset{\varminimizer}{\min} \ \underset{\varmaximizer}{\max}\ \hat{\mathcal{L}}(\varminimizer, \varmaximizer, \auxvarminimizeritr{k}, \auxvarmaximizeritr{k}, \dualvarminimizeritr{k}, \dualvarmaximizeritr{k}) \) \label{step:saddleADMMvarupdate}
        
         \(\auxvarminimizeritr{k+1} = \arg \underset{\auxvarminimizer}{\min} \ \hat{\mathcal{L}}(\varminimizeritr{k+1}, \varmaximizeritr{k+1}, \auxvarminimizer, \auxvarmaximizeritr{k}, \dualvarminimizeritr{k}, \dualvarmaximizeritr{k})\) \label{step:saddleADMMauxvarupdate1}
         
        \(\auxvarmaximizeritr{k+1} = \arg \underset{\auxvarmaximizer}{\max} \ \hat{\mathcal{L}}(\varminimizeritr{k+1}, \varmaximizeritr{k+1}, \auxvarminimizeritr{k+1}, \auxvarmaximizer, \dualvarminimizeritr{k}, \dualvarmaximizeritr{k})\) \label{step:saddleADMMauxvarupdate2}
         
        \(\dualvarminimizeritr{k+1} = \dualvarminimizeritr{k} + \dualstepsizeminimizer(\varminimizeritr{k+1} - \auxvarminimizeritr{k+1})\)
        \(\dualvarmaximizeritr{k+1} = \dualvarmaximizeritr{k} + \dualstepsizemaximizer(\varmaximizeritr{k+1} - \auxvarmaximizeritr{k+1})\) \label{step:saddleADMMdualupdate}

    }

    \caption{Saddle-Point Alternating Direction Method of Multipliers (SP-ADMM) for decomposable optimization}
    \label{alg:sp-admm}
\end{algorithm}

Saddle-point ADMM for decomposable optimization, Algorithm \ref{alg:sp-admm}, also consists of three steps: primal variable \(\varminimizer, \varmaximizer\)  updates, auxiliary primal variable \(\auxvarminimizer, \auxvarmaximizer\) updates, and dual variable \(\dualvarminimizer,\dualvarmaximizer\) updates. Line \ref{step:saddleADMMvarupdate} of Algorithm \ref{alg:sp-admm} is separable: This step assigns 
\begin{align*}
    \arg &\min_{\varminimizerblock{i} \in \spaceminimizerblock{i}} \ \max_{\varmaximizerblock{i} \in \spacemaximizerblock{i}}  \funcblock{i}(\varminimizerblock{i}, \varmaximizerblock{i})
    \\
    &+(\dualvarminimizerblockitr{i}{k})^{\top}(\varminimizerblock{i} - \auxvarminimizerblockitr{i}{k}) + \frac{\dualstepsizeminimizer}{2} \ltwonorm{\varminimizerblock{i} - \auxvarminimizerblockitr{i}{k}}
    \\
    &- (\dualvarmaximizerblockitr{i}{k})^{\top}(\varmaximizerblock{i} - \auxvarmaximizerblockitr{i}{k}) - \frac{\dualstepsizemaximizer}{2} \ltwonorm{\varmaximizerblock{i} - \auxvarmaximizerblockitr{i}{k}}
\end{align*} to \(\varminimizerblockitr{i}{k+1}\) and \(\varmaximizerblockitr{i}{k+1}\) for every \(i \in 1, \ldots N\). These sub-problems have a significantly lower number of dimensions compared to the original saddle-point problem \eqref{prob:original} and can be solved in parallel. The sub-problems  can be solved using existing saddle-point optimization methods and for some objective functions such as bilinear functions of two one-dimensional variables, they have analytical solutions. Lines \ref{step:saddleADMMauxvarupdate1}--\ref{step:saddleADMMauxvarupdate2} are the convex projection steps and are equal to letting \[\auxvarminimizeritr{k+1} = \arg \min_{\auxvarminimizer \in \spaceminimizer} \ltwonorm{\varminimizeritr{k+1} + \dualvarminimizeritr{k}/\dualstepsizeminimizer - \auxvarminimizer} \] and \[\auxvarmaximizeritr{k+1} = \arg \min_{\auxvarmaximizer \in \spacemaximizer} \ltwonorm{\varmaximizeritr{k+1} + \dualvarmaximizeritr{k}/\dualstepsizemaximizer - \auxvarmaximizer}, \] which can be solved using convex optimization methods. 

We show the convergence of SP-ADMM, under a similar assumption of standard ADMM. We assume that there exists a saddle-point where strong duality holds for the minimizer's problem when the maximizer is fixed, and vice versa.
\begin{assumption} \label{ass:saddlepoint}
There exists \((\varminimizeritr{*}, \varmaximizeritr{*}, \auxvarminimizeritr{*}, \auxvarmaximizeritr{*}, \dualvarminimizeritr{*}, \dualvarmaximizeritr{*})\) such that 

\begin{align}
    &\mathcal{L} (\varminimizeritr{*}, \varmaximizeritr{*}, \auxvarminimizeritr{*}, \auxvarmaximizeritr{*}, \dualvarminimizer, \dualvarmaximizeritr{*}) \nonumber
    \\
    &\leq \mathcal{L} (\varminimizeritr{*}, \varmaximizeritr{*}, \auxvarminimizeritr{*}, \auxvarmaximizeritr{*}, \dualvarminimizeritr{*}, \dualvarmaximizeritr{*}) \label{ineq:saddlepointmin}
    \\
    &\leq     \mathcal{L} (\varminimizer, \varmaximizeritr{*}, \auxvarminimizer, \auxvarmaximizeritr{*}, \dualvarminimizeritr{*}, \dualvarmaximizeritr{*}) \nonumber
\end{align}
and
\begin{align} 
    &\mathcal{L} (\varminimizeritr{*}, \varmaximizer, \auxvarminimizeritr{*}, \auxvarmaximizer, \dualvarminimizeritr{*}, \dualvarmaximizeritr{*}) \nonumber
    \\
    &\leq \mathcal{L} (\varminimizeritr{*}, \varmaximizeritr{*}, \auxvarminimizeritr{*}, \auxvarmaximizeritr{*}, \dualvarminimizeritr{*}, \dualvarmaximizeritr{*})  \label{ineq:saddlepointmax}
    \\
    &\leq     \mathcal{L} (\varminimizeritr{*}, \varmaximizeritr{*}, \auxvarminimizeritr{*}, \auxvarmaximizeritr{*}, \dualvarminimizeritr{*}, \dualvarmaximizer) \nonumber
\end{align}
for all \(\varminimizer, \varmaximizer, \auxvarminimizer, \auxvarmaximizer, \dualvarminimizer,\) and \(\dualvarmaximizer\).
\end{assumption} 
Note that \(\varminimizeritr{*} = \auxvarminimizeritr{*}\) and \(\varmaximizeritr{*} = \auxvarmaximizeritr{*}\) for the saddle-point since \(\sup_{\dualvarminimizer} \dualvarminimizer^{\top}(\varminimizeritr{*} - \auxvarminimizeritr{*}) = \infty\) and \(\inf_{\dualvarmaximizer} \dualvarmaximizer^{\top}(\varmaximizeritr{*} - \auxvarmaximizeritr{*})= -\infty\) otherwise. Also note that \(\varminimizeritr{*} \in \spaceminimizer \cap (\spaceminimizerblock{1}\times \ldots \times \spaceminimizerblock{N})\) and \(\varmaximizeritr{*} \in \spacemaximizer \cap (\spacemaximizerblock{1}\times \ldots \times \spacemaximizerblock{N})\) due to the indicator functions.

Despite its complicated nature, the assumption is satisfied for the convex-concave saddle-point point problems where Slater's condition~\cite{boyd2004convex} is satisfied. 
\begin{proposition}[Sufficient condition for a saddle-point] \label{prop:existenceofsaddlepoint}
There exists a saddle point \((\varminimizeritr{*}, \varmaximizeritr{*}, \auxvarminimizeritr{*}, \auxvarmaximizeritr{*}, \dualvarminimizeritr{*}, \dualvarmaximizeritr{*})\) for \(\mathcal{L}\) that satisfies Assumption \ref{ass:saddlepoint} if 
\begin{enumerate}
    \item Every \(\funcblock{i}\) is a convex function of \(\varminimizerblock{i}\) and concave function of \(\varmaximizerblock{i}\) in \(\spaceminimizerblock{i}\times \spacemaximizerblock{i}\).     
    \label{condition:convexconcaveclosedobjective}
    \item Every \(\funcblock{i}\) is continuous. \label{condition:continuity}
    \item \(\spaceminimizer\), \(\spacemaximizer\), and every \(\spaceminimizerblock{i}\), \(\spacemaximizerblock{i}\) are compact, convex polytopes. \label{condition:compactconstraints}
\end{enumerate}
\end{proposition}

We give the proof of Proposition \ref{prop:existenceofsaddlepoint} in Appendix \ref{apx:existenceofsaddlepoint}.

Note that the conditions given in Proposition \ref{prop:existenceofsaddlepoint} imply that the saddle-point problem satisfies Slater's condition for the minimizer and maximizer. In the proposition, we use polytope constraints for simplicity; the proposition can be improved to general convex sets  \(\spaceminimizerblock{i}\), \(\spacemaximizerblock{i}\), \(\spaceminimizer\), and \(\spacemaximizer\) as long as there is a saddle point for \eqref{prob:original} that satisfies Slater's condition. 

Under Assumption \ref{ass:saddlepoint}, the iterates of SP-ADMM converges to a saddle point of \eqref{prob:original}. If every \(\funcblock{i}\) is Lipschitz continuous, the proposition also implies the convergence of value. 
\begin{proposition} \label{prop:convergence}
 Under Assumption \ref{ass:saddlepoint}, the iterates of SP-ADMM converge to a saddle point for \eqref{prob:original}, i.e., \(\auxvarminimizeritr{k} \to \varminimizeritr{*}\) and \(\auxvarmaximizeritr{k} \to \varmaximizeritr{*}\)  where \((\varminimizeritr{*}, \varmaximizeritr{*})\) is a saddle-point of \eqref{prob:original}.
\end{proposition}
We give the proof of Proposition \ref{prop:convergence} in Appendix \ref{apx:convergence}.

Proposition \ref{prop:convergence} shows convergence in the limit. As in the standard ADMM~\cite{boyd2011distributed}, one can use the magnitude of primal and dual residuals as the stopping criterion in practice: Terminate when \( \ltwonormnormal{\varminimizeritr{k} - \auxvarminimizeritr{k}} + \ltwonormnormal{\varmaximizeritr{k} - \auxvarmaximizeritr{k}} \leq \epsilon^{\text{primal}}\) and \(  \dualstepsizeminimizer \ltwonormnormal{\auxvarminimizeritr{k} - \auxvarminimizeritr{k-1}} + \dualstepsizemaximizer \ltwonormnormal{\auxvarmaximizeritr{k} - \auxvarmaximizeritr{k-1}} \leq \epsilon^{\text{dual}}\) where \(\varminimizeritr{k} - \auxvarminimizeritr{k}\) and \(\varmaximizeritr{k} - \auxvarmaximizeritr{k}\) are the primal residuals, and \(\dualstepsizeminimizer (\auxvarminimizeritr{k} - \auxvarminimizeritr{k-1})\) and \( \dualstepsizemaximizer(\auxvarmaximizeritr{k} - \auxvarmaximizeritr{k-1})\) are the dual residuals after iteration \(k\).

\section{Numerical Examples}
In this section, we give numerical examples for SP-ADMM and compare it with saddle-point Frank-Wolfe (SP-FW) method~\cite{gidel2017frank}. The implementations are given at \url{https://github.com/mustafakarabag/SP-ADMM}.
\subsection{Power Allocation Game for Communication Channels}

\begin{figure}[t]
       \centering

    \input{power_allocation_game_val}
    \caption{(Top) Total capacity of the communication channels with \((\auxvarminimizeritr{k}, \auxvarmaximizeritr{k})\). (Bottom) Total residual norm is \(\ltwonormnormal{\varminimizeritr{k} - \auxvarminimizeritr{k}} + \ltwonormnormal{\varmaximizeritr{k} - \auxvarmaximizeritr{k}} + \dualstepsizeminimizer \ltwonormnormal{\auxvarminimizeritr{k} - \auxvarminimizeritr{k-1}} + \dualstepsizemaximizer \ltwonormnormal{\auxvarmaximizeritr{k} - \auxvarmaximizeritr{k-1}}.\)}
    \label{fig:power_allocation_example}
\end{figure}

In this example from \cite{ghosh2003minimax}, we consider a power allocation problem in Gaussian communication channels. The total communication capacity is \(\sum_{i=1}^{N} \log\left(1 + \frac{\varmaximizerblock{i}}{\sigma_{i} + \varminimizerblock{i}}\right)\) where \(\varmaximizerblock{i}\) is the signal power allocated to the \(i^{\text{th}}\) channel, \(\sigma_{i}\) is the receiver noise for the \(i^{\text{th}}\) channel, and \(\varminimizerblock{i}\) is the  noise of the \(i^{\text{th}}\) channel. 

We consider a game between a maximizer that allocates signal powers and a minimizer that adversarially chooses the noise levels for \(N=10\) channels. The global constraints are  \(\sum_{i=1}^{N} \varmaximizerblock{i} = 20\) for the the maximizer and \(\sum_{i=1}^{N} \varminimizerblock{i} = 10\) for the minimizer. Players have individual constraints \(\varminimizerblock{i} \geq 0\) and \(\varmaximizerblock{i} \geq 0\). The receiver noise level is \(\sigma = [2, 6, 5, 8, 3, 9, 5, 6, 7, 3]\). The equilibrium value of the problem instance is \(2.860\)~\cite{ghosh2003minimax}.

For the implementation of SP-ADMM, we use SP-FW to solve the sub-saddle-point problems that are in the form of 
\begin{align*}
    &\min_{\varminimizerblock{i}} \max_{\varmaximizerblock{i}} \sum_{i=1}^{N} \log\left(1 + \frac{\varmaximizerblock{i}}{\sigma_{i} + \varminimizerblock{i}}\right) + \dualvarminimizerblockitr{i}{k}(\varminimizerblock{i} - \auxvarminimizerblockitr{i}{k}) 
    \\
    &+ \dualstepsizeminimizer(\varminimizerblock{i} - \auxvarminimizerblockitr{i}{k})^2  - \dualvarmaximizerblockitr{i}{k}(\varmaximizerblock{i} - \auxvarmaximizerblockitr{i}{k}) - \dualstepsizemaximizer(\varmaximizerblock{i} - \auxvarmaximizerblockitr{i}{k})^2.
\end{align*} We initialize \(\varminimizer\) and \(\varmaximizer\) with a vector of zeros. The variables \(\auxvarminimizer\) and \(\auxvarmaximizer\) are initialized with the projections of \(\varminimizer\) and \(\varmaximizer\) onto their global constraints, respectively.

In Figure \ref{fig:power_allocation_example}, we show the output of SP-ADMM for different penalty parameters. Similar to the standard ADMM, SP-ADMM generates acceptable solutions within a few iterations: The total capacity converges to the equilibrium value 2.860. The total residual norm decay as the number of iterations increase. However, similar to the standard ADMM, the rate of convergence is slow. We suspect that the fluctuations of the total residual norm is due to the dynamic competition between the players and the fact that sub-problems are solved with a finite accuracy. When we compare the effects of the penalty parameters \(\dualstepsizeminimizer\) and \(\dualstepsizemaximizer\), we observe that mild penalties such as \(0.1\) lead to both faster objective and residual convergences.

\subsection{Network Routing Game with Adversarial Agents}
In this example, we consider a network routing problem represented with a Markov decision process (MDP). The MDP is deterministic, i.e., it is a directed graph with \(N\) edges. Players choose a policy for this MDP that induces a Markov chain. The players' policies control the density of atomic agents that are transitioning in the Markov chain. The variables, \(\varminimizer\) and \(\varmaximizer\), of the players represent the stationary distributions induced by the players over the edges of the Markov chain. We generate the underlying directed graph of the MDP using a random Erdos-Renyi graph such that every node has 5 edges in expectation.

The network has a price function for every edge \(i\) that is equal to \(\varminimizerblock{i} + \varmaximizerblock{i}\), i.e., the total demand for edge \(i\). The cost of an edge \(i\), for the minimizer is \(\varminimizerblock{i}(\varminimizerblock{i} + \varmaximizerblock{i})\) that is the density of minimizer times the price of the edge. The minimizer's goal is to minimize the total cost \(\sum_{i=1}^{N} \varminimizerblock{i}(\varminimizerblock{i} + \varmaximizerblock{i})\). The maximizer is an adversary trying to maximize the same cost. The minimizer and maximizer control a unit density each. The individual constraints are \(0 \leq \varminimizerblock{i} \leq 1\) and \(0 \leq \varmaximizerblock{i} \leq 1\) for every edge \(i\).  The global contraints are enforced by the dynamics of the MDP: The players' stationary distributions have to be valid. In addition, the maximizer's density at state \(1\) has to be at least 0.1, i.e., \(\sum_{i \in E} \varminimizerblock{i} \geq 0.1\) where \(E\) is the incoming edges of state \(1\).

We compare the performance of SP-ADMM with SP-FW for different sizes of MDPs. For the initialization of both SP-ADMM with SP-FW, we use the valid stationary distribution that is closest to the uniform distribution in \(L_{2}\) distance. We solve the sub-saddle-point problems of SP-ADMM using an analytical solution exploiting the bilinear structure of sub-problems. This step has \(\mathcal{O}(N)\) time complexity. The gradients for SP-FW are also computed using analytical solutions, which has \(\mathcal{O}(N)\) time complexity. The projection step of SP-ADMM and the maximization step of SP-FW are both computed using ECOS solver~\cite{domahidi2013ecos} with CVXPY~\cite{diamond2016cvxpy} interface. For SP-ADMM, we use \(\dualstepsizeminimizer=\dualstepsizemaximizer=1\), and for SP-FW, we use the step size \(2/(2+k)\) at iteration \(k\) as suggested in \cite{gidel2017frank}.

For both algorithms, we compute a bound on the optimality gap in the following way. Let \(\auxvarminimizeritr{*,k}\) be the optimal response of the minimizer against the maximizer's \(\auxvarmaximizeritr{k}\) action, and \(\auxvarmaximizeritr{*,k}\) be the optimal response of the maximizer against the minimizer's \(\auxvarminimizeritr{k}\) action. We compute the best action of a player by solving a convex optimization problem where the other player's action is fixed. By the definition of a saddle-point, we have \[\sum_{i=1}^{N} \funcblock{i}(\auxvarminimizerblockitr{i}{*,k}, \auxvarmaximizerblockitr{i}{k}) \leq \sum_{i=1}^{N} \funcblock{i}(\auxvarminimizerblockitr{i}{*}, \auxvarmaximizerblockitr{i}{*}) \leq \sum_{i=1}^{N} \funcblock{i}(\auxvarminimizerblockitr{i}{k}, \auxvarmaximizerblockitr{i}{*,k}). \]
The best lower bound is \(l^{k} = \max_{1 \leq j \leq k} \sum_{i=1}^{N} \funcblock{i}(\auxvarminimizerblockitr{i}{*,j}, \auxvarmaximizerblockitr{i}{j}) \) and best upper bound is \(u^{k} = \min_{1 \leq j \leq k} \sum_{i=1}^{N} \funcblock{i}(\auxvarminimizerblockitr{i}{j}, \auxvarmaximizerblockitr{i}{*,j})\) at iteration \(k\). The optimality gap of an iterative algorithm at iteration \(k\) is bounded by \(u^{k} - l^{k}\). 

\begin{table}[t]
\caption{Comparison of SP-ADMM and SP-FW for different MDP sizes}
\begin{center}
\centering
\begin{tabular}{|cc|cc|cc|}
\hline
\multicolumn{2}{|c|}{Network size}                                                                       & \multicolumn{2}{c|}{SP-ADMM}                                                                     & \multicolumn{2}{c|}{SP-FW}                                                                       \\ \hline
\multicolumn{1}{|c|}{\# nodes} & \begin{tabular}[c]{@{}c@{}} \# edges \\ \(N\)\end{tabular} & \begin{tabular}[c]{@{}c@{}}Opt. gap\\ \(u^{k} - l^{k}\)\end{tabular} & Time (s) & \begin{tabular}[c]{@{}c@{}}Opt. gap\\ \(u^{k} - l^{k}\)\end{tabular} & Time (s) \\ \hline
\multicolumn{1}{|c|}{10}              & 49                                                               & 1.36e-9$^{\mathrm{a}}$                                                                       & 5.48             & 1.36e-9$^{\mathrm{a}}$                                                                    & 5.17             \\
\multicolumn{1}{|c|}{20}              & 93                                                               & 2.49e-7                                                                       & 9.39             & 5.13e-3                                                                      & 9.09            \\
\multicolumn{1}{|c|}{50}              & 282                                                              & 1.87e-6                                                                       & 28.06            & 2.35e-3                                                                       & 25.67            \\
\multicolumn{1}{|c|}{100}             & 494                                                              & 1.35e-6                                                                       & 51.18           & 1.62e-3                                                                       & 48.17            \\ \hline
\multicolumn{6}{l}{
\begin{tabular}{@{}l@{}}$^{\mathrm{a}}$ Both algorithms fail to improve on the initialization point\\   due to numerical precision issues.\end{tabular}
}
\end{tabular}

\label{tab:sp-admmvssp-fw}
\end{center}
\end{table}

\begin{figure}[t]
    \centering
    \input{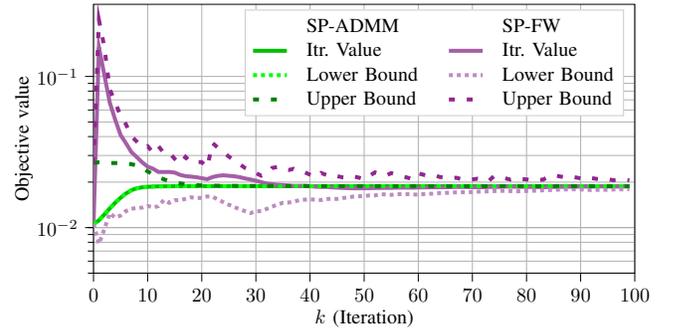}
    \caption{The objective values for SP-ADMM and SP-FW. For each algorithm, 'Itr. Value' refers to the value with the variables from the current iterate, i.e., \((\auxvarminimizeritr{k}, \auxvarmaximizeritr{k,*})\). 'Lower Bound' refers to the value with the maximizer's variable from the current iterate and the minimizer's best response to it, i.e., \((\auxvarminimizeritr{k,*}, \auxvarmaximizeritr{k})\). 'Upper Bound' refers to the value with the minimizer's variable from the current iterate and the maximizer's best response to it, i.e., \((\auxvarminimizeritr{k}, \auxvarmaximizeritr{k,*})\). }
    \label{fig:sp-admmvssp-fw}
\end{figure}

We compare SP-ADMM and SP-FW in Table \ref{tab:sp-admmvssp-fw} and Figure \ref{fig:sp-admmvssp-fw}. In Figure \ref{fig:sp-admmvssp-fw}, we observe that SP-ADMM performs better than SP-FW for objective convergence. In addition, the upper and lower bounds are closer for SP-ADMM, which shows a better convergence to the saddle-point solution. In Table \ref{tab:sp-admmvssp-fw}, we observe that the solution time for SP-ADMM is slightly worse since we solve a quadratic program for SP-ADMM whereas we solve a linear program of the same size for SP-FW. On the other hand, the optimality gap \(u^{k}-l^{k}\) is orders of magnitude better for SP-ADMM with similar solution times.

\section{Conclusion}
We demonstrated saddle-point alternating direction method of multipliers (SP-ADMM) to solve decomposable saddle-point problems. We show that SP-ADMM has convergence guarantees under a saddle-point assumption. This assumption is satisfied for convex-concave problems that satisfy Slater's conditions. While we show that SP-ADMM converges asymptotically, we suspect that it also enjoys the non-asymptotic guarantees of standard ADMM~\cite{nishihara2015general}, for example, in the strongly convex-strongly concave setting.

\bibliographystyle{IEEEtran}
\bibliography{ref}

\appendices

\section{Proof of Proposition \ref{prop:existenceofsaddlepoint}} \label{apx:existenceofsaddlepoint}
We show the existence of a saddle point for the augmented Lagrangian by considering the minimax theorem~\cite{border1985fixed} and Slater's constraint qualification for convex duality~\cite{boyd2004convex}.
Since \(\sum_{i=1}^{N} \funcblock{i}(\varminimizerblock{i}, \varmaximizerblock{i})\) is a continous, convex-concave function and the feasible spaces are compact, convex for both minimizer and maximizer, there exists a saddle point \((\varminimizeritr{*}, \varmaximizeritr{*})\) for \eqref{prob:original} by the minimax theorem~\cite{border1985fixed}. Consequently, \((\varminimizeritr{*}, \varmaximizeritr{*}, \auxvarminimizeritr{*}, \auxvarmaximizeritr{*})\) is a saddle point of \eqref{prob:withaux} where \(\auxvarminimizeritr{*} = \varminimizeritr{*}\) and \(\varmaximizeritr{*} = \auxvarmaximizeritr{*}\).

Since the feasible spaces are compact, convex polytopes, 
\begin{itemize}
    \item there exist \(G_{\minimizer,i}\) and \(h_{\minimizer,i}\) such that \(G_{\minimizer,i}\varminimizerblock{i} + h_{\minimizer,i} \leq 0\) is equal to \(\varminimizerblock{i} \in \spaceminimizerblock{i}\),
    \item there exist \(G_{\maximizer,i}\) and \(h_{\maximizer,i}\) such that \(G_{\maximizer,i}\varmaximizerblock{i} + h_{\maximizer,i} \leq 0\) is equal to \(\varmaximizerblock{i} \in \spacemaximizerblock{i}\),
\item there exist \(G_{\minimizer}\) and \(h_{\minimizer}\) such that \(G_{\minimizer}\auxvarminimizer + h_{\minimizer} \leq 0\) is equal to \(\auxvarminimizer \in \spaceminimizer\), and
    \item there exist \(G_{\maximizer}\) and \(h_{\maximizer}\) such that \(G_{\maximizer}\auxvarmaximizer + h_{\maximizer} \leq 0\) is equal to \(\auxvarmaximizer \in \spacemaximizer.\)

\end{itemize}

Define the Lagrangian for \eqref{prob:withaux}
\begin{align*} 
    & \bar{\mathcal{L}}(\varminimizer, \varmaximizer, \auxvarminimizer, \auxvarmaximizer, \dualvarminimizer, \dualvarmaximizer, \mu_{\minimizer}, \mu_{\maximizer}, [ \mu_{\minimizer, i}]_{i=1}^{N},  [ \mu_{\maximizer, i}]_{i=1}^{N}) 
    \\
    =&  \sum_{i=1}^{N} \funcblock{i}(\varminimizerblock{i}, \varmaximizerblock{i}) 
    \\
    +&  \dualvarminimizer^{\top}(\varminimizer - \auxvarminimizer)  +  \mu_{\minimizer}^{\top} (G_{\minimizer}\auxvarminimizer + h_{\minimizer}) + \sum_{i=1}^{N}  \mu_{\minimizer, i}^{\top} (G_{\minimizer,i}\varminimizerblock{i} + h_{\minimizer,i}) 
    \\
    -&  \dualvarmaximizer^{\top}(\varmaximizer - \auxvarmaximizer) - \mu_{\maximizer}^{\top} (G_{\maximizer}\auxvarmaximizer + h_{\maximizer}) - \sum_{i=1}^{N} \mu_{\maximizer, i}^{\top} (G_{\maximizer,i}\varmaximizerblock{i} + h_{\maximizer,i})
\end{align*}
where \(\mu_{\minimizer}, \mu_{\maximizer}, \mu_{\minimizer, 1}, \mydots, \mu_{\minimizer, N}, \mu_{\maximizer, 1}, \mydots, \mu_{\maximizer, N} \geq 0\).

For fixed \(\varmaximizeritr{*}\) and \(\auxvarmaximizeritr{*}\),  \(\sum_{i=1}^{N} \funcblock{i}(\varminimizerblock{i}, \varmaximizerblockitr{i}{*})\) is a continuous, jointly convex function of \(\varminimizer\) and \(\dualvarminimizer\) and the constraints of \eqref{prob:withaux} satisfies Slater's condition. Note that \(\varminimizeritr{*}\) and \(\auxvarminimizeritr{*}\) is optimal for fixed \(\varmaximizeritr{*}\) and \(\auxvarmaximizeritr{*}\). By the saddle point theorem~\cite{burke}, there exists \( (\varminimizeritr{*}, \auxvarminimizeritr{*} , \dualvarminimizeritr{*}, \mu_{\minimizer}^{*}, \mu_{\minimizer, 1}^{*}, \mydots, \mu_{\minimizer, N}^{*}) \) such that
\begin{subequations} \label{saddlebarlagrangianminimizer}
\begin{align}
        &\bar{\mathcal{L}}(\varminimizeritr{*}, \varmaximizeritr{*}, \auxvarminimizeritr{*}, \auxvarmaximizeritr{*}, \dualvarminimizer, \dualvarmaximizer, \mu_{\minimizer}, \mu_{\maximizer}, [ \mu_{\minimizer, i}]_{i=1}^{N}, [ \mu_{\maximizer, i}]_{i=1}^{N})  \label{saddlebarlagrangianminimizer1}
        \\
        &\leq\bar{\mathcal{L}}(\varminimizeritr{*}, \varmaximizeritr{*}, \auxvarminimizeritr{*}, \auxvarmaximizeritr{*}, \dualvarminimizeritr{*}, \dualvarmaximizer, \mu_{\minimizer}^{*}, \mu_{\maximizer}, [ \mu_{\minimizer, i}^{*}]_{i=1}^{N}, [ \mu_{\maximizer, i}]_{i=1}^{N})  \label{saddlebarlagrangianminimizer2}
        \\
        &\leq\bar{\mathcal{L}}(\varminimizeritr{*}, \varmaximizeritr{*}, \auxvarminimizeritr{*}, \auxvarmaximizeritr{*}, \dualvarminimizeritr{*}, \dualvarmaximizer, \mu_{\minimizer}^{*}, \mu_{\maximizer}, [ \mu_{\minimizer, i}^{*}]_{i=1}^{N}, [ \mu_{\maximizer, i}]_{i=1}^{N})  \label{saddlebarlagrangianminimizer3}
\end{align}
\end{subequations} for any \(\dualvarmaximizer, \mu_{\maximizer},  [ \mu_{\maximizer, i}]_{i=1}^{N}\). Let \(\dualvarmaximizer, \mu_{\maximizer},  [ \mu_{\maximizer, i}]_{i=1}^{N} = 0\). 
Note that 
\begin{align*}
     &\mathcal{L} (\varminimizeritr{*}, \varmaximizeritr{*}, \auxvarminimizeritr{*}, \auxvarmaximizeritr{*}, \dualvarminimizeritr{*}, \dualvarmaximizer) 
     \\
     &\leq \bar{\mathcal{L}}(\varminimizeritr{*}, \varmaximizeritr{*}, \auxvarminimizeritr{*}, \auxvarmaximizeritr{*}, \dualvarminimizeritr{*}, \dualvarmaximizer, \mu_{\minimizer}, 0, [ \mu_{\minimizer, i}]_{i=1}^{N},[ 0]_{i=1}^{N}) 
\end{align*} since \(G_{\minimizer,i}\varminimizerblockitr{i}{*} + h_{\minimizer,i} \leq 0\), \(G_{\minimizer,i}\auxvarminimizeritr{*} + h_{\minimizer,i} \leq 0\), and \(\mu_{\minimizer}, \mu_{\minimizer, 1}, \mydots, \mu_{\minimizer, N} \geq 0\). We also have  \begin{align*}
     &\bar{\mathcal{L}}(\varminimizer, \varmaximizeritr{*}, \auxvarminimizer, \auxvarmaximizeritr{*}, \dualvarminimizeritr{*}, \dualvarmaximizer, \mu_{\minimizer}^{*}, 0, [ \mu_{\minimizer, i}^{*}]_{i=1}^{N}, [ 0]_{i=1}^{N}) 
     \\
     &\leq \mathcal{L} (\varminimizer, \varmaximizeritr{*}, \auxvarminimizer, \auxvarmaximizeritr{*}, \dualvarminimizeritr{*}, \dualvarmaximizer)
\end{align*} since \(\indicatorfunc{\spaceminimizerblock{i}}  (\varminimizerblock{i}) \geq  \mu_{\minimizer, i}^{\top} (G_{\minimizer,i}\varminimizerblock{i} + h_{\minimizer,i})\), \(\indicatorfunc{\spaceminimizer}  (\auxvarminimizer) \geq  \mu_{\minimizer}^{\top} (G_{\minimizer}\varminimizer + h_{\minimizer})\), \(\varmaximizerblockitr{i}{*} \in \spacemaximizerblock{i}\), and \(\auxvarmaximizeritr{*} \in \spacemaximizer\).

We established \[\mathcal{L} (\varminimizeritr{*}, \varmaximizeritr{*}, \auxvarminimizeritr{*}, \auxvarmaximizeritr{*}, \dualvarminimizeritr{*}, \dualvarmaximizer) \leq \mathcal{L} (\varminimizer, \varmaximizeritr{*}, \auxvarminimizer, \auxvarmaximizeritr{*}, \dualvarminimizeritr{*}, \dualvarmaximizer).\] We now show \[\mathcal{L}(\varminimizeritr{*}, \varmaximizeritr{*}, \auxvarminimizeritr{*}, \auxvarmaximizeritr{*}, \dualvarminimizer,\dualvarmaximizer) \leq \mathcal{L}(\varminimizeritr{*}, \varmaximizeritr{*}, \auxvarminimizeritr{*}, \auxvarmaximizeritr{*}, \dualvarminimizeritr{*},\dualvarmaximizer).\] Note that the optimization problem \[\min_{\dualvarminimizer, \mu_{\minimizer},  [ \mu_{\minimizer, i}]_{i=1}^{N}}\bar{\mathcal{L}}(\varminimizeritr{*}, \varmaximizeritr{*}, \auxvarminimizeritr{*}, \auxvarmaximizeritr{*}, \dualvarminimizer, \dualvarmaximizer, \mu_{\minimizer}, \mu_{\maximizer}, [ \mu_{\minimizer, i}]_{i=1}^{N}, [ \mu_{\maximizer, i}]_{i=1}^{N}) \] is separable: the optimal values of \(\dualvarminimizer\) and \(\mu_{\minimizer},  \mu_{\minimizer, 1}, \mydots, \mu_{\minimizer, N}\) can be computed independently. Consequently, since \(\dualvarminimizer^{*}\) is a maximizer for \[\bar{\mathcal{L}}(\varminimizeritr{*}, \varmaximizeritr{*}, \auxvarminimizeritr{*}, \auxvarmaximizeritr{*}, \dualvarminimizer, \dualvarmaximizer, \mu_{\minimizer}, \mu_{\maximizer}, [ \mu_{\minimizer, i}]_{i=1}^{N}, [ \mu_{\maximizer, i}]_{i=1}^{N}),\] it is also a maximizer for \(\mathcal{L}(\varminimizeritr{*}, \varmaximizeritr{*}, \auxvarminimizeritr{*}, \auxvarmaximizeritr{*}, \dualvarminimizer,\dualvarmaximizer),\) and we have \(\mathcal{L}(\varminimizeritr{*}, \varmaximizeritr{*}, \auxvarminimizeritr{*}, \auxvarmaximizeritr{*}, \dualvarminimizer,\dualvarmaximizer) \leq \mathcal{L}(\varminimizeritr{*}, \varmaximizeritr{*}, \auxvarminimizeritr{*}, \auxvarmaximizeritr{*}, \dualvarminimizeritr{*},\dualvarmaximizer).\) Combining these results, we get 
\begin{align}
    &\mathcal{L} (\varminimizeritr{*}, \varmaximizeritr{*}, \auxvarminimizeritr{*}, \auxvarmaximizeritr{*}, \dualvarminimizer, \dualvarmaximizer) \nonumber
    \\
    &\leq \mathcal{L} (\varminimizeritr{*}, \varmaximizeritr{*}, \auxvarminimizeritr{*}, \auxvarmaximizeritr{*}, \dualvarminimizeritr{*}, \dualvarmaximizer) \label{saddlebarlagrangianminimizerfreemaximizer}
    \\
    &\leq     \mathcal{L} (\varminimizer, \varmaximizeritr{*}, \auxvarminimizer, \auxvarmaximizeritr{*}, \dualvarminimizeritr{*}, \dualvarmaximizer) \nonumber
\end{align}
 for arbitrary \(\dualvarmaximizer\). By symmetry, we can repeat the same arguments and get 
\begin{align}
    &\mathcal{L} (\varminimizeritr{*}, \varmaximizer, \auxvarminimizeritr{*}, \auxvarmaximizer, \dualvarminimizer, \dualvarmaximizeritr{*}) \nonumber
    \\
    &\leq \mathcal{L} (\varminimizeritr{*}, \varmaximizeritr{*}, \auxvarminimizeritr{*}, \auxvarmaximizeritr{*}, \dualvarminimizer, \dualvarmaximizeritr{*})  \label{saddlebarlagrangianmaximizerfreeminimizer}
    \\
    &\leq     \mathcal{L} (\varminimizeritr{*}, \varmaximizeritr{*}, \auxvarminimizeritr{*}, \auxvarmaximizeritr{*}, \dualvarminimizer, \dualvarmaximizer) \nonumber
\end{align}
for arbitrary \(\dualvarminimizer\).
Finally, by letting \(\dualvarmaximizer = \dualvarmaximizeritr{*}\) in \eqref{saddlebarlagrangianminimizerfreemaximizer} and \(\dualvarminimizer = \dualvarminimizeritr{*}\) in \eqref{saddlebarlagrangianmaximizerfreeminimizer}, we get the desired result.

\section{Proof of Proposition \ref{prop:convergence}} \label{apx:convergence}
The proof follows the same steps of the proof for convergence for the standard ADMM algorithm~\cite{boyd2011distributed}. The work \cite{boyd2011distributed} proves convergence of standard ADMM by considering only the properties of minimizer updates. To prove the convergence of SP-ADMM, we consider the properties of both minimizer and maximizer updates.

We define the value function of the algorithm  \[\valuefunctionitr{k} = \frac{\ltwonorm{\dualvarminimizeritr{k} - \dualvarminimizeritr{*}}}{\dualstepsizeminimizer} + \frac{\ltwonorm{\dualvarmaximizeritr{k} - \dualvarmaximizeritr{*}}}{\dualstepsizemaximizer} + \frac{\ltwonorm{\auxvarminimizeritr{k} - \auxvarminimizeritr{*}}}{1/\dualstepsizeminimizer} +  \frac{\ltwonorm{\auxvarmaximizeritr{k} - \auxvarmaximizeritr{*}}}{1/\dualstepsizemaximizer}.\]  We will show that the value decreases at every step, i.e.,  \begin{align} 
    \valuefunctionitr{k+1} \leq& \valuefunctionitr{k} - \dualstepsizeminimizer \ltwonorm{\primalresidualminimizeritr{k+1}} - \dualstepsizemaximizer \ltwonorm{\primalresidualmaximizeritr{k+1}} \nonumber
    \\ 
    &- \dualstepsizeminimizer \ltwonorm{\auxvarminimizeritr{k+1} - \auxvarminimizeritr{k}} - \dualstepsizemaximizer \ltwonorm{\auxvarmaximizeritr{k+1} - \auxvarmaximizeritr{k}}. \label{ineq:valuedecreases}
\end{align} where \(\primalresidualminimizeritr{k}= \varminimizeritr{k} - \auxvarminimizeritr{k}\) is the primal residual for the minimizer and \(\primalresidualmaximizeritr{k}= \varmaximizeritr{k} - \auxvarmaximizeritr{k}\) is the primal residual for the maximizer. By telescoping sum over \(k\), we get 
\begin{align} 
    \valuefunctionitr{0} \geq& \sum_{i=1}^{\infty} \dualstepsizeminimizer \ltwonorm{\primalresidualminimizeritr{k}} + \dualstepsizemaximizer \ltwonorm{\primalresidualmaximizeritr{k}} \nonumber
    \\
    &+ \dualstepsizeminimizer \ltwonorm{\auxvarminimizeritr{k} - \auxvarminimizeritr{k-1}} + \dualstepsizemaximizer \ltwonorm{\auxvarmaximizeritr{k} - \auxvarmaximizeritr{k-1}}. \nonumber
\end{align}

Since \(\valuefunctionitr{0}\) is finite, and \(\dualstepsizeminimizer\) and \(\dualstepsizemaximizer\) are strictly positive, we must have \(\lim_{k \to \infty} \ltwonorm{\primalresidualminimizeritr{k}} = 0 \),  \(\lim_{k \to \infty} \ltwonorm{\primalresidualmaximizeritr{k}} = 0 \),  \(\lim_{k \to \infty} \ltwonorm{\auxvarminimizeritr{k} - \auxvarminimizeritr{k-1}} = 0 \), and  \(\lim_{k \to \infty} \ltwonorm{\auxvarmaximizeritr{k} - \auxvarmaximizeritr{k-1}} = 0 \). Consequently, \(\varminimizeritr{k} \to \varminimizeritr{\circ}\), \(\auxvarminimizeritr{k} \to \varminimizeritr{\circ}\), \(\varmaximizeritr{k} \to \varmaximizeritr{\circ}\), and \(\auxvarmaximizeritr{k} \to \varmaximizeritr{\circ}\) for some stationary point \((\varminimizeritr{\circ},\varmaximizeritr{\circ},\varminimizeritr{\circ},\varmaximizeritr{\circ},\dualvarminimizeritr{\circ},\dualvarmaximizeritr{\circ})\). Since \((\varminimizeritr{\circ},\varmaximizeritr{\circ},\varminimizeritr{\circ},\varmaximizeritr{\circ},\dualvarminimizeritr{\circ},\dualvarmaximizeritr{\circ})\) is a stationary point of SP-ADMM, \((\varminimizeritr{\circ},\varminimizeritr{\circ},\dualvarminimizeritr{\circ})\) is a stationary point of ADMM (Algorithm \ref{alg:admm}) when \((\varmaximizer=\varmaximizeritr{\circ},\auxvarmaximizer=\varmaximizeritr{\circ},\dualvarmaximizer=\dualvarmaximizeritr{\circ})\) is fixed, and therefore \((\varminimizeritr{\circ}, \varminimizeritr{\circ})\) is a solution to \eqref{prob:withaux} when \(\varmaximizer=\varmaximizeritr{\circ},\auxvarmaximizer=\varmaximizeritr{\circ}\) is fixed. Similarly, \((\varmaximizeritr{\circ}, \varmaximizeritr{\circ})\) is a solution to \eqref{prob:withaux} when \(\varminimizer=\varminimizeritr{\circ},\auxvarminimizer=\varminimizeritr{\circ}\) is fixed. Consequently, \((\varminimizeritr{\circ},\varmaximizeritr{\circ},\varminimizeritr{\circ},\varmaximizeritr{\circ})\) is a saddle point of \eqref{prob:withaux} and \((\varminimizeritr{\circ},\varmaximizeritr{\circ})\) is a saddle point of \eqref{prob:original}.

We now show \eqref{ineq:valuedecreases}. For ease of notation, we also define the following quantities:
\begin{itemize}  
                   \item  Equilibrium value \(\optvalue = \sum_{i=1}^{N} \funcblock{i}(\varminimizerblockitr{i}{*}, \varmaximizerblockitr{i}{*})\). Note that \(\optvalue = \mathcal{L}(\varminimizeritr{*}, \varmaximizeritr{*}, \auxvarminimizeritr{*}, \auxvarmaximizeritr{*}, \dualvarminimizeritr{*}, \dualvarmaximizeritr{*}) \) since \(\varminimizeritr{*} = \auxvarminimizeritr{*},\) \(\varmaximizeritr{*} = \auxvarmaximizeritr{*}\)
                \item \(\valueitr{k} = \sum_{i=1}^{N} \funcblock{i}(\varminimizerblockitr{i}{k}, \varmaximizerblockitr{i}{k})\), \(\optvalueminimizeritr{k} = \sum_{i=1}^{N} \funcblock{i}(\varminimizerblockitr{i}{k}, \varmaximizerblockitr{i}{*})\),  \(\optvaluemaximizeritr{k} = \sum_{i=1}^{N} \funcblock{i}(\varminimizerblockitr{i}{*}, \varmaximizerblockitr{i}{k})\).
\end{itemize}

To prove \eqref{ineq:valuedecreases}, we will show 
\begin{equation} \label{ineq:minimizerconvergesfrombottom}
    \optvalue - \optvalueminimizeritr{k+1} \leq (\dualvarminimizeritr{*})^{\top} \primalresidualminimizeritr{k+1},
\end{equation}
\begin{equation} \label{ineq:maximizerconvergesfrombottom}
     \optvaluemaximizeritr{k+1} - \optvalue \leq (\dualvarmaximizeritr{*})^{\top} \primalresidualmaximizeritr{k+1},
\end{equation}
\begin{align} 
        \valueitr{k+1} - \optvaluemaximizeritr{k+1}   \leq &  \dualstepsizeminimizer (\auxvarminimizeritr{k+1} - \auxvarminimizeritr{k})^{\top} (\primalresidualminimizeritr{k+1} + \auxvarminimizeritr{k+1} - \auxvarminimizeritr{*})  \nonumber
        \\
        &  -(\dualvarminimizeritr{k+1})^{\top} \primalresidualminimizeritr{k+1}, \label{ineq:minimizerconvergesfromtop}
\end{align}
and 
\begin{align} 
        \optvalueminimizeritr{k+1} - \valueitr{k+1}   \leq& \dualstepsizemaximizer (\auxvarmaximizeritr{k+1} - \auxvarmaximizeritr{k})^{\top} (\primalresidualmaximizeritr{k+1} + \auxvarmaximizeritr{k+1} - \auxvarmaximizeritr{*}) \nonumber
        \\
        &-(\dualvarmaximizeritr{k+1})^{\top} \primalresidualmaximizeritr{k+1}. \label{ineq:maximizerconvergesfromtop}
\end{align} We, for now, assume that these inequalities hold and give the proofs in 
Appendix \ref{proof:ineqset1} and \ref{proof:ineqset2}.

\subsection{Proof of \eqref{ineq:valuedecreases}}

Adding \eqref{ineq:minimizerconvergesfrombottom}, \eqref{ineq:maximizerconvergesfrombottom}, \eqref{ineq:minimizerconvergesfromtop}, and \eqref{ineq:maximizerconvergesfromtop}, and multiplying by 2, we get 
\begin{align}
    0 \leq & 2(\dualvarminimizeritr{*}-\dualvarminimizeritr{k+1})^{\top} \primalresidualminimizeritr{k+1} + 2(\dualvarmaximizeritr{*}-\dualvarmaximizeritr{k+1})^{\top} \primalresidualmaximizeritr{k+1} \nonumber
    \\
    &+ 2\dualstepsizeminimizer (\auxvarminimizeritr{k+1} - \auxvarminimizeritr{k})^{\top} (\primalresidualminimizeritr{k+1} + \auxvarminimizeritr{k+1} - \auxvarminimizeritr{*}) \nonumber
    \\
    &+ 2\dualstepsizemaximizer (\auxvarmaximizeritr{k+1} - \auxvarmaximizeritr{k})^{\top} (\primalresidualmaximizeritr{k+1} + \auxvarmaximizeritr{k+1} - \auxvarmaximizeritr{*}). \label{ineq:4addedterms}
\end{align} We use the definitions to rewrite \eqref{ineq:4addedterms}.

Using \(\dualvarminimizeritr{k+1} = \dualvarminimizeritr{k} + \dualstepsizeminimizer \primalresidualminimizeritr{k+1}\), \(\primalresidualminimizeritr{k+1} = (\dualvarminimizeritr{k+1} - \dualvarminimizeritr{k})/\dualstepsizeminimizer\), \(\dualvarminimizeritr{k+1} -\dualvarminimizeritr{k}=\dualvarminimizeritr{k+1}-\dualvarminimizeritr{*} + \dualvarminimizeritr{*}-\dualvarminimizeritr{k} \), we get
{\small
\begin{align}
    2&(\dualvarminimizeritr{*}-\dualvarminimizeritr{k+1})^{\top} \primalresidualminimizeritr{k+1} = 2(\dualvarminimizeritr{*}-\dualvarminimizeritr{k})^{\top} \primalresidualminimizeritr{k+1} - 2\dualstepsizeminimizer \ltwonorm{\primalresidualminimizeritr{k+1}} \nonumber
    \\
    =&\frac{2}{\dualstepsizeminimizer}(\dualvarminimizeritr{*} - \dualvarminimizeritr{k})^{\top} (\dualvarminimizeritr{k+1} - \dualvarminimizeritr{*}) - \frac{1}{\dualstepsizeminimizer} \ltwonorm{\dualvarminimizeritr{k+1} - \dualvarminimizeritr{k}}- \dualstepsizeminimizer \ltwonorm{\primalresidualminimizeritr{k+1}} \nonumber
    \\
    =&\frac{1}{\dualstepsizeminimizer}\ltwonorm{\dualvarminimizeritr{k} - \dualvarminimizeritr{*}} - \frac{1}{\dualstepsizeminimizer}\ltwonorm{\dualvarminimizeritr{k+1} - \dualvarminimizeritr{*}} - \dualstepsizeminimizer \ltwonorm{\primalresidualminimizeritr{k+1}}. \label{eq:replacement1stterm}
\end{align}}%
By the symmetry of the definitions, we also get
\begin{align}
        &2(\dualvarmaximizeritr{*}-\dualvarmaximizeritr{k+1})^{\top} \primalresidualmaximizeritr{k+1} 
    \\
    &=\frac{1}{\dualstepsizemaximizer}\ltwonorm{\dualvarmaximizeritr{k} - \dualvarmaximizeritr{*}} - \frac{1}{\dualstepsizemaximizer}\ltwonorm{\dualvarmaximizeritr{k+1} - \dualvarmaximizeritr{*}} - \dualstepsizemaximizer \ltwonorm{\primalresidualmaximizeritr{k+1}} \label{eq:replacement2ndterm}
\end{align}

Using \(\auxvarminimizeritr{k+1} - \auxvarminimizeritr{*} = \auxvarminimizeritr{k+1} - \auxvarminimizeritr{k} + \auxvarminimizeritr{k}- \auxvarminimizeritr{*}\) and \(\auxvarminimizeritr{k+1} - \auxvarminimizeritr{k} = \auxvarminimizeritr{k+1} - \auxvarminimizeritr{*} - \auxvarminimizeritr{k}+ \auxvarminimizeritr{*}\), we get
\begin{align}
    2&\dualstepsizeminimizer (\auxvarminimizeritr{k+1} - \auxvarminimizeritr{k})^{\top} (\primalresidualminimizeritr{k+1} + \auxvarminimizeritr{k+1} - \auxvarminimizeritr{*}) - \ltwonorm{\primalresidualminimizeritr{k+1}} \nonumber
    \\
    =&-\dualstepsizeminimizer\ltwonorm{\primalresidualminimizeritr{k+1} + \auxvarminimizeritr{k+1} - \auxvarminimizeritr{k}} \nonumber
    \\
    &-\dualstepsizeminimizer(\ltwonorm{ \auxvarminimizeritr{k+1} - \auxvarminimizeritr{*}} - \ltwonorm{ \auxvarminimizeritr{k} - \auxvarminimizeritr{*}}) \label{eq:replacement3rdterm}
\end{align}
By the symmetry of the definitions, we also get
\begin{align}
    2&\dualstepsizemaximizer (\auxvarmaximizeritr{k+1} - \auxvarmaximizeritr{k})^{\top} (\primalresidualmaximizeritr{k+1} + \auxvarmaximizeritr{k+1} - \auxvarmaximizeritr{*}) - \ltwonorm{\primalresidualmaximizeritr{k+1}} \nonumber
    \\
    =&-\dualstepsizemaximizer\ltwonorm{\primalresidualmaximizeritr{k+1} + \auxvarmaximizeritr{k+1} - \auxvarmaximizeritr{k}} \nonumber
    \\
    &-\dualstepsizemaximizer(\ltwonorm{ \auxvarmaximizeritr{k+1} - \auxvarmaximizeritr{*}} - \ltwonorm{ \auxvarmaximizeritr{k} - \auxvarmaximizeritr{*}}). \label{eq:replacement4thterm}
\end{align}
By substituting \eqref{eq:replacement1stterm}, \eqref{eq:replacement2ndterm}, \eqref{eq:replacement3rdterm}, and \eqref{eq:replacement4thterm} in \eqref{ineq:4addedterms}, we get 
\begin{align}
    &\valuefunctionitr{k+1} \leq \valuefunctionitr{k} - \dualstepsizeminimizer\ltwonorm{\primalresidualminimizeritr{k+1} + \auxvarminimizeritr{k+1} - \auxvarminimizeritr{k}} 
    \\
    &\quad \quad \quad -\dualstepsizemaximizer\ltwonorm{\primalresidualmaximizeritr{k+1} + \auxvarmaximizeritr{k+1} - \auxvarmaximizeritr{k}} \nonumber
    \\
    &\leq \valuefunctionitr{k} - \dualstepsizeminimizer \ltwonorm{\primalresidualminimizeritr{k+1}} - \dualstepsizemaximizer \ltwonorm{\primalresidualmaximizeritr{k+1}}  \nonumber
    \\
    &- \dualstepsizeminimizer \ltwonorm{\auxvarminimizeritr{k+1} - \auxvarminimizeritr{k}} - \dualstepsizemaximizer \ltwonorm{\auxvarmaximizeritr{k+1} - \auxvarmaximizeritr{k}}\nonumber
    \\
    &-2\dualstepsizeminimizer(\dualvarminimizeritr{k+1})^{\top}(\auxvarminimizeritr{k+1} - \auxvarminimizeritr{k}) -2\dualstepsizemaximizer(\dualvarmaximizeritr{k+1})^{\top}(\auxvarmaximizeritr{k+1} - \auxvarmaximizeritr{k}) \label{ineq:valuefunctionlastform}
\end{align}

As shown in the proofs of \eqref{ineq:minimizerconvergesfromtop} and \eqref{ineq:maximizerconvergesfromtop}, \(\auxvarminimizeritr{k+1}\) minimizes \(-(\dualvarminimizeritr{k+1})^{\top}\auxvarminimizer\) in \(\spaceminimizer\), and \(\auxvarmaximizeritr{k+1}\) maximizes \((\dualvarmaximizeritr{k+1})^{\top}\auxvarmaximizer\) in \(\spacemaximizer\). Consequently, we have \(-2\dualstepsizeminimizer(\dualvarminimizeritr{k+1})^{\top}\auxvarminimizeritr{k+1} \leq -2\dualstepsizeminimizer(\dualvarminimizeritr{k+1})^{\top} \auxvarminimizeritr{k},\) and \(2\dualstepsizemaximizer(\dualvarmaximizeritr{k+1})^{\top}\auxvarmaximizeritr{k+1} \leq -2\dualstepsizemaximizer(\dualvarmaximizeritr{k+1})^{\top} \auxvarmaximizeritr{k}.\) Combining these with \eqref{ineq:valuefunctionlastform}, we get \eqref{ineq:valuedecreases}.

\subsection{Proofs of \eqref{ineq:minimizerconvergesfrombottom} and \eqref{ineq:maximizerconvergesfrombottom}}   \label{proof:ineqset1}
Due to the saddle point assumption, we have \[\mathcal{L}(\varminimizeritr{*}, \varmaximizeritr{*}, \auxvarminimizeritr{*}, \auxvarmaximizeritr{*}, \dualvarminimizeritr{*}, \dualvarmaximizeritr{*}) \leq \mathcal{L}(\varminimizeritr{k+1}, \varmaximizeritr{*}, \auxvarminimizeritr{k+1}, \auxvarmaximizeritr{*}, \dualvarminimizeritr{*}, \dualvarmaximizeritr{*}).\] Since \(\optvalue = \mathcal{L} (\varminimizeritr{*}, \varmaximizeritr{*}, \auxvarminimizeritr{*}, \auxvarmaximizeritr{*}, \dualvarminimizeritr{*}, \dualvarmaximizeritr{*})\) and \(\varmaximizeritr{*} = \auxvarmaximizeritr{*}\), we have \[\optvalue \leq \optvalueminimizeritr{k+1} + (\dualvarminimizeritr{*})^{\top} (\varminimizeritr{k+1} - \auxvarminimizeritr{k+1}).\] Using \(\primalresidualminimizeritr{k+1} = \varminimizeritr{k+1} - \auxvarminimizeritr{k+1}\) and rearranging the terms, we get \begin{equation} \label{ineq:minimizerconvergesfrombottomproof}
      \optvalue - \optvalueminimizeritr{k+1} \leq (\dualvarminimizeritr{*})^{\top} \primalresidualminimizeritr{k+1}.
\end{equation}  

The proof of \eqref{ineq:maximizerconvergesfrombottom} has the same steps with the proof of \eqref{ineq:minimizerconvergesfrombottom}. 

\subsection{Proofs of \eqref{ineq:minimizerconvergesfromtop} and \eqref{ineq:maximizerconvergesfromtop}}
\label{proof:ineqset2}

We note that \(\hat{\mathcal{L}}(\varminimizer, \varmaximizer, \auxvarminimizeritr{k}, \auxvarmaximizeritr{k}, \dualvarminimizeritr{k}, \dualvarmaximizeritr{k})\) is a convex function of \(\varminimizer\) and a concave function of \(\varmaximizer\), and \((\varminimizeritr{k+1}, \varmaximizeritr{k+1})\) is a solution to \[ \min_{\varminimizer \in \spaceminimizerblock{1} \times \ldots \times \spaceminimizerblock{N}} \ \max_{\varmaximizer \in \spacemaximizerblock{1} \times \ldots \times \spacemaximizerblock{N}} \ \hat{\mathcal{L}}(\varminimizer, \varmaximizer, \auxvarminimizeritr{k}, \auxvarmaximizeritr{k}, \dualvarminimizeritr{k}, \dualvarmaximizeritr{k}).\]

Define 
\begin{align*}
    g(\varminimizer, \varmaximizer) =&   \sum_{i=1}^{N} \funcblock{i}(\varminimizerblock{i}, \varmaximizerblock{i}) +   (\dualvarminimizeritr{k} - \dualstepsizeminimizer (\auxvarminimizeritr{k+1} - \auxvarminimizeritr{k}))^{\top} \varminimizer 
    \\
    &- (\dualvarmaximizeritr{k} - \dualstepsizemaximizer (\auxvarmaximizeritr{k+1} - \auxvarmaximizeritr{k}))^{\top} \varmaximizer
\end{align*}

Using \(\dualvarminimizeritr{k+1} \in \spaceminimizerblock{1} \times \ldots \times \spaceminimizerblock{N}\) and \(\dualvarminimizeritr{k+1} = \dualvarminimizeritr{k} + \dualstepsizeminimizer (\varminimizeritr{k+1} - \auxvarminimizeritr{k+1})\), we get \begin{align*}
    \frac{\partial \hat{\mathcal{L}}(\varminimizer, \varmaximizer, \auxvarminimizeritr{k}, \auxvarmaximizeritr{k}, \dualvarminimizeritr{k}, \dualvarmaximizeritr{k})}{ \partial \varminimizer}\Biggr|_{\substack{\varminimizer = \varminimizeritr{k+1}}} =\frac{\partial g(\varminimizer, \varmaximizer)}{ \partial \varminimizer}\Biggr|_{\substack{\varminimizer = \varminimizeritr{k+1}}}
\end{align*}
Similarly, we get \[\frac{\partial \hat{\mathcal{L}}(\varminimizer, \varmaximizer, \auxvarminimizeritr{k}, \auxvarmaximizeritr{k}, \dualvarminimizeritr{k}, \dualvarmaximizeritr{k})}{ \partial \varmaximizer}\Biggr|_{\substack{\varmaximizer = \varmaximizeritr{k+1}}} = \frac{\partial g(\varminimizer, \varmaximizer)}{ \partial \varmaximizer}\Biggr|_{\substack{\varmaximizer = \varmaximizeritr{k+1}}}.\] Since \(\hat{\mathcal{L}}(\varminimizer, \varmaximizer, \auxvarminimizeritr{k}, \auxvarmaximizeritr{k}, \dualvarminimizeritr{k}, \dualvarmaximizeritr{k})\) and \( g(\varminimizer, \varmaximizer)\) share the same gradient field for \(\varminimizer\) and \(\varmaximizer\), and \((\varminimizeritr{k+1}, \varmaximizeritr{k+1})\) is a saddle point of \(\hat{\mathcal{L}}(\varminimizer, \varmaximizer, \auxvarminimizeritr{k}, \auxvarmaximizeritr{k}, \dualvarminimizeritr{k}, \dualvarmaximizeritr{k})\),  \((\varminimizeritr{k+1}, \varmaximizeritr{k+1})\) is also a saddle point of \( g(\varminimizer, \varmaximizer)\). Using the saddle point property we have,
\begin{align*}
    &\sum_{i=1}^{N} \funcblock{i}(\varminimizerblockitr{i}{k+1}, \varmaximizerblockitr{i}{k+1}) +   (\dualvarminimizeritr{k+1} - \dualstepsizeminimizer (\auxvarminimizeritr{k+1} - \auxvarminimizeritr{k}))^{\top} \varminimizeritr{k+1} 
    \\
    & \quad \quad- (\dualvarmaximizeritr{k+1} - \dualstepsizemaximizer (\auxvarmaximizeritr{k+1} - \auxvarmaximizeritr{k}))^{\top} \varmaximizeritr{k+1} 
    \\
    &\leq \sum_{i=1}^{N} \funcblock{i}(\varminimizerblockitr{i}{*}, \varmaximizerblockitr{i}{k+1})+ (\dualvarminimizeritr{k+1} - \dualstepsizeminimizer (\auxvarminimizeritr{k+1} - \auxvarminimizeritr{k}))^{\top} \varminimizeritr{*}  
    \\
    & \quad \quad- (\dualvarmaximizeritr{k+1} - \dualstepsizemaximizer (\auxvarmaximizeritr{k+1} - \auxvarmaximizeritr{k}))^{\top} \varmaximizeritr{k+1}
\end{align*}
By definitions of \(\optvaluemaximizeritr{k+1}\) and \(\valueitr{k+1}\), we get 
\begin{align} 
    \valueitr{k+1} +   &(\dualvarminimizeritr{k+1} - \dualstepsizeminimizer (\auxvarminimizeritr{k+1} - \auxvarminimizeritr{k}))^{\top} \varminimizeritr{k+1} \nonumber
    \\
    &\leq \optvaluemaximizeritr{k+1} + (\dualvarminimizeritr{k+1} - \dualstepsizeminimizer (\auxvarminimizeritr{k+1} - \auxvarminimizeritr{k}))^{\top} \varminimizeritr{*}. \label{ineq:tobecombined1}
\end{align}
By the saddle point property, we also get 
\begin{align}
    \valueitr{k+1} - &(\dualvarmaximizeritr{k+1} - \dualstepsizemaximizer (\auxvarmaximizeritr{k+1} - \auxvarmaximizeritr{k}))^{\top} \varmaximizeritr{k+1} \nonumber
    \\
    &\geq \optvalueminimizeritr{k+1} - (\dualvarmaximizeritr{k+1} - \dualstepsizemaximizer (\auxvarmaximizeritr{k+1} - \auxvarmaximizeritr{k}))^{\top} \varmaximizeritr{*}. \label{ineq:tobecombined2}
\end{align}

Define \(h_{\minimizer}(\auxvarminimizer) =  -(\dualvarminimizeritr{k+1})^{\top}\auxvarminimizer.\) and \(h_{\maximizer}(\auxvarmaximizer) =  (\dualvarmaximizeritr{k+1})^{\top}\auxvarmaximizer.\) We have
\begin{align*}
    \frac{\partial \hat{\mathcal{L}}(\varminimizer^{k+1}, \varmaximizer^{k+1}, \auxvarminimizer, \auxvarmaximizeritr{k}, \dualvarminimizeritr{k}, \dualvarmaximizeritr{k})}{ \partial \auxvarminimizer}\Biggr|_{\substack{\auxvarminimizer = \auxvarminimizeritr{k+1}}}  =\frac{\partial  h_{\minimizer}(\auxvarminimizer) }{\partial \varminimizer}\Biggr|_{\substack{\varminimizer = \auxvarminimizeritr{k+1}}}
\end{align*} and similarly 
\begin{align*}
    \frac{\partial \hat{\mathcal{L}}(\varminimizer^{k+1}, \varmaximizer^{k+1}, \auxvarminimizeritr{k+1}, \auxvarmaximizer, \dualvarminimizeritr{k}, \dualvarmaximizeritr{k})}{ \partial \auxvarmaximizer}\Biggr|_{\substack{\auxvarmaximizer = \auxvarmaximizeritr{k+1}}} & =\frac{\partial  h_{\maximizer}(\auxvarmaximizer) }{\partial \varmaximizer}\Biggr|_{\substack{\varmaximizer = \auxvarmaximizeritr{k+1}}}.
\end{align*}
Since \(\hat{\mathcal{L}}(\varminimizer^{k+1}, \varmaximizer^{k+1}, \auxvarminimizer, \auxvarmaximizeritr{k}, \dualvarminimizeritr{k}, \dualvarmaximizeritr{k})\) and \(h_{\minimizer}(\auxvarminimizer)\) has the same gradient field in \(\spaceminimizer\), \(\auxvarminimizeritr{k+1}\) is also a minimizer of \(h_{\minimizer}(\auxvarminimizer)\) in \(\spaceminimizer\). Similarly, \(\auxvarmaximizeritr{k+1}\) is also a maximizer of \(h_{\maximizer}(\auxvarmaximizer)\) in \(\spacemaximizer\). Due to these we have
\begin{equation} \label{ineq:tobecombined3}
    -(\dualvarminimizeritr{k+1})^{\top}\auxvarminimizeritr{k+1} \leq -(\dualvarminimizeritr{k+1})^{\top}\auxvarminimizeritr{*}
\end{equation}
and 
\begin{equation} \label{ineq:tobecombined4}
    (\dualvarmaximizeritr{k+1})^{\top}\auxvarmaximizeritr{k+1} \geq (\dualvarmaximizeritr{k+1})^{\top}\auxvarmaximizeritr{*}.
\end{equation}
By combining \eqref{ineq:tobecombined1} and \eqref{ineq:tobecombined3}, and noting that \(\varminimizeritr{*} = \auxvarminimizeritr{*}\) and \(\primalresidualminimizeritr{k+1} = \varminimizeritr{k+1} - \auxvarminimizeritr{k+1}\), we get \begin{align} 
        \valueitr{k+1} - \optvaluemaximizeritr{k+1}   \leq &  \dualstepsizeminimizer (\auxvarminimizeritr{k+1} - \auxvarminimizeritr{k})^{\top} (\primalresidualminimizeritr{k+1} + \auxvarminimizeritr{k+1} - \auxvarminimizeritr{*})  \nonumber
        \\
        &  -(\dualvarminimizeritr{k+1})^{\top} \primalresidualminimizeritr{k+1} \label{ineq:minimizerconvergesfromtopproof}
\end{align}
Similarly, by combining \eqref{ineq:tobecombined2} and \eqref{ineq:tobecombined4}, and noting that \(\varmaximizeritr{*} = \auxvarmaximizeritr{*}\) and \(\primalresidualmaximizeritr{k+1} = \varmaximizeritr{k+1} - \auxvarmaximizeritr{k+1}\), we get \begin{align} 
        \optvalueminimizeritr{k+1} - \valueitr{k+1}   \leq& \dualstepsizemaximizer (\auxvarmaximizeritr{k+1} - \auxvarmaximizeritr{k})^{\top} (\primalresidualmaximizeritr{k+1} + \auxvarmaximizeritr{k+1} - \auxvarmaximizeritr{*}) \nonumber
        \\
        &-(\dualvarmaximizeritr{k+1})^{\top} \primalresidualmaximizeritr{k+1}. \label{ineq:maximizerconvergesfromtopproof}
\end{align}

\end{document}